\documentclass[12pt]{article}
\usepackage{amssymb}
\usepackage{amsxtra}
\usepackage{amstext}
\usepackage{latexsym}
\usepackage{amsmath}
\usepackage{dsfont}
\usepackage{mathrsfs}

\DeclareMathAlphabet{\mathpzc}{OT1}{pzc}{m}{it}

\newtheorem{definition}{Definition}[section]
\newtheorem{theorem}[definition]{Theorem}
\newtheorem{lemma}[definition]{Lemma}

\newtheorem{note}[definition]{Note}
\newtheorem{assumption}[definition]{Assumption}
\newtheorem{notation}[definition]{Notation}

\def\K{\mathbb K}

\begin{document}

\title{\bf Bidiagonal Triads and the Tetrahedron Algebra}
\author{Darren Funk-Neubauer}
\date{}

\maketitle

\begin{abstract}
\noindent
We introduce a linear algebraic object called a bidiagonal triad.  A bidiagonal triad is a modification of the previously studied and similarly defined concept of bidiagonal triple.  A bidiagonal triad and a bidiagonal triple both consist of three diagonalizable linear transformations on a finite-dimensional vector space, each of which acts in a bidiagonal fashion on the eigenspaces of the other two.  A triad differs from a triple in the way these bidiagonal actions are defined.  We modify a number of theorems about bidiagonal triples to the case of  bidiagonal triads.  We also describe the close relationship between bidiagonal triads and the representation theory of the tetrahedron Lie algebra.            \\


\end{abstract}

\section{Introduction}

In this paper we introduce a linear algebraic object called a bidiagonal triad.  Roughly speaking, a bidiagonal triad is an ordered list of three diagonalizable linear transformations on a finite-dimensional vector space, each of which acts in a bidiagonal fashion on the eigenspaces of the other two (see Definition \ref{def:bdtriad} for the precise definition).  This concept arose as a modification of the similarly defined concept of a bidiagonal triple (see Definition \ref{def:bdtriple}).  Furthermore, the idea of a bidiagonal triple arose as a generalization of the concept of a bidiagonal pair (see Definition \ref{def:bdpair}).  The theory of bidiagonal pairs and triples was developed in \cite{Funk-Neubauer13, Funk-Neubauer17}.  The relationship between bidiagonal pairs, triples, and triads is as follows.  If $A,\, A', \, A''$ is a bidiagonal triple then $A,\, A'$ and $A',\, A''$ and $A'',\, A$ are each bidiagonal pairs.  However, the converse is not true.  That is, if $A,\, A'$ and $A',\, A''$ and $A'',\, A$ are each bidiagonal pairs then $A,\, A', \, A''$ does {\it not} have to be a bidiagonal triple.  This is the genesis of the concept of bidiagonal triad.  Three bidiagonal pairs can be interwoven together to produce a bidiagonal triple, but if the interweaving process is altered then a bidiagonal triad will be produced.  For a triad the three transformations have a common eigenvector.  However, for a triple there does not exist
a common eigenvector for all three transformations, but any two of the three have a common eigenvector.   \\

The main purpose of this paper is twofold.  First we show how some results about bidiagonal triples in \cite{Funk-Neubauer17} can be modified to hold for bidiagonal triads.  In particular, we show that the eigenvalues of the three linear transformations in a bidiagonal triad satisfy a linear recurrence (see Theorem \ref{thm:recur}).  We then define a type of bidiagonal triad called reduced based on the solutions to this recurrence (see Definition \ref{def:reduced}).  We show that reduced bidiagonal triads are canonical in the sense that every bidiagonal triad is equivalent to a reduced one (see Definition \ref{def:ae} and Theorem \ref{thm:reduce}).  The second main purpose of this paper is to explain the connection between bidiagonal triads and representations of the tetrahedron algebra.  In Theorem \ref{thm:tetBDtriad} we show that every finite-dimensional irreducible module for the tetrahedron algebra provides numerous examples of reduced bidiagonal triads.  Our last main result is a partial converse of Theorem \ref{thm:tetBDtriad}.  In Theorem \ref{thm:BDtriadtet} we show how an irreducible module for the tetrahedron algebra can be constructed starting from a certain type of reduced bidiagonal triad called thin (see Definition \ref{def:thin}).  The stronger non-thin version of Theorem \ref{thm:BDtriadtet} is false as we demonstrate with a counterexample at the end of the paper.    \\

We now offer some background on the tetrahedron algebra.  The tetrahedron algebra, denoted $\boxtimes$, was introduced in \cite{HarTer07}.  $\boxtimes$ has essentially six generators, each of which can be identified with the six edges of a tetrahedron.  In other words, each generator is indexed by a pair of vertices in the tetrahedron.  The symmetric group $S_4$ acts on $\boxtimes$ as a group of automorphisms \cite[Section 2]{HarTer07}.  $\boxtimes$ is isomorphic to the three-point $\mathfrak{sl}_2$ loop algebra \cite[Theorem 11.5]{HarTer07}.  Any five of the six edges of the tetrahedron generate a subalgebra of $\boxtimes$ which is isomorphic to the $\mathfrak{sl}_2$ loop algebra \cite[Corollary 12.6]{HarTer07}.  The three surrounding edges in each of the four faces of the tetrahedron form a basis for a subalgebra of $\boxtimes$ isomorphic to $\mathfrak{sl}_2$ \cite[Corollary 12.4]{HarTer07}.  Each of the three pairs of opposite edges of the tetrahedron generate a subalgebra of $\boxtimes$ isomorphic to the Onsager algebra \cite[Corollary 12.5]{HarTer07}.  $\boxtimes$ decomposes into the direct sum of these three Onsager subalgebras \cite[Theorem 11.6]{HarTer07}.  For more information on the algebras mentioned in the previous five sentences see the references in \cite[Section 1]{HarTer07}.  The finite-dimensional irreducible $\boxtimes$-modules were classified in \cite[Theorems 1.7, 1.8]{Hartwig07}.  This classification invokes the classification of the finite-dimensional irreducible modules for the Onsager algebra (see \cite{DateRoan, Davies90, Davies91} as well as the summary in \cite[Theorems 1.3, 1.4, 1.6]{Hartwig07}).  For more information on how the tetrahedron algebra arises in representation theory see \cite{Al-Najjar13, Benkart07, Elduque07, ItoTerinpress4, NeherSavSen}.  See \cite{HouGao, Morales, MorPas} for how the tetrahedron algebra has been used in algebraic combinatorics and graph theory.         \\

We now offer some additional information on bidiagonal pairs and triples because of their close connection to bidiagonal triads.  Bidiagonal pairs and triples have been used to study representations of the well-known algebras $\mathfrak{sl}_2$, $U_q(\mathfrak{sl}_2)$, $U_q(\widehat{\mathfrak{sl}}_2)$, $\boxtimes$, and $\boxtimes_q$.  The $q$-tetrahedron algebra $\boxtimes_q$ is a quantum analog of the tetrahedron algebra \cite{ItoTer072}.  Bidiagonal pairs and triples appear implicitly in \cite{Benkart04, Funk-Neubauer07, Funk-Neubauer09, ItoTer073, ItoTerWang06, Ter13}.  However, they were not explicitly defined until \cite{Funk-Neubauer13, Funk-Neubauer17} which contain systematic studies of bidiagonal pairs and triples as well as classification theorems.  Finite-dimensional irreducible $\boxtimes$-modules provide numerous examples of bidiagonal pairs and triples.  The actions of the three generators corresponding to the surrounding edges of each face of the tetrahedron act as a bidiagonal triple \cite[Theorem 4.8]{Funk-Neubauer17}, \cite[Theorem 3.8]{Hartwig07}, \cite[Corollary 12.4]{HarTer07}.  The actions of the two generators corresponding to any two edges of the tetrahedron that meet in a common vertex act as a bidiagonal pair \cite[Theorem 5.8]{Funk-Neubauer17}.  See also Note \ref{BDtriplestriadspairs} of this paper.  Similarly, finite-dimensional irreducible $\boxtimes_q$-modules also provide numerous examples of bidiagonal pairs and triples \cite{Funk-Neubauer17, ItoTer072}.  See \cite{Al-Najjar20, Baseilhac21, Lynch, Nomura21, Ter14, Ter15, Ter17, Yang16, Yang20} for more references to bidiagonal pairs and triples.        \\

We now remark on the initial motivation for investigating bidiagonal pairs, triples, and triads as well as the tetrahedron algebra.  These objects each had their genesis in the study of tridiagonal pairs.  Roughly speaking, a tridiagonal pair is an ordered pair of diagonalizable linear transformations on a finite-dimensional vector space, each of which acts in a tridiagonal fashion on the eigenspaces of the other (see \cite[Definition 1.1]{ItoTanTer01} for the precise definition).  Tridiagonal pairs originally arose in algebraic combinatorics through the study of a combinatorial object called a P- and Q-polynomial association scheme \cite{ItoTanTer01}.  Tridiagonal pairs appear in a wide variety of mathematical subjects including representation theory, orthogonal polynomials, special functions, partially ordered sets, statistical mechanics, and classical mechanics.  See \cite[Introduction]{Funk-Neubauer13} for the appropriate references.  A major classification result involving tridiagonal pairs appears in \cite{ItoTer09}.  Finite-dimensional irreducible $\boxtimes$-modules provide numerous examples of tridiagonal pairs.  The actions of each pair of generators corresponding to opposite edges of the tetrahedron act as a tridiagonal pair \cite[Theorem 1.7, Corollary 2.7]{Hartwig07}.  These tridiagonal pairs are said to have Krawtchouk type \cite{ItoTer074}.  See also Note \ref{BDtriplestriadspairs} of this paper.  Similarly, finite-dimensional irreducible $\boxtimes_q$-modules also provide numerous examples of tridiagonal pairs \cite[Theorem 10.3]{ItoTer072}, \cite[Theorem 2.7]{ItoTer075}.  These tridiagonal pairs are said to have $q$-geometric type \cite[Definition 2.6]{ItoTer075}.  A certain special case of a tridiagonal pair, called a Leonard pair, has also been the subject of much research.  The Leonard pairs are classified up to isomorphism \cite[Theorem 1.9]{Ter011}.  See \cite{Ter042} for a survey of Leonard pairs.  \\

This paper is organized as follows.  In Section 2 we recall the definitions of bidiagonal pair and triple and then define bidiagonal triad along with a number of related definitions.  In Section 3 we recall the definition of the tetrahedron algebra $\boxtimes$ and review some facts about $\boxtimes$-modules.  We also describe the relationship between $\boxtimes$ and the Lie algebra $\mathfrak{sl}_2$.  Section 4 contains statements of the main results of the paper.  In Section 5 we develop the tools needed to prove our main results including some additional properties of bidiagonal triads.  Sections 6 through 10 contain the proofs of the main results stated in Section 4.  In Section 11 we provide a counterexample to show that a stronger version of Theorem \ref{thm:BDtriadtet} is false.  

\section{Bidiagonal Triads}

\noindent
\begin{notation}
\label{note1}
\rm
Throughout this paper we adopt the following notation.  Let $d$ denote a nonnegative integer.  Let $\K$ denote a field.  Let $V$ denote a vector space over $\K$ with finite positive dimension.  For linear transformations $X: V \to V$ and $Y: V \to V$ we define $[X, Y] := XY - YX$.
\end{notation}

\begin{note}
\label{brackprops}
\rm
With reference to Notation \ref{note1} observe that $[X,X] = 0$ and $[X,Y] = -[Y,X]$.  Also, $[\,\, , \,\,]$ is bilinear.  That is, for $c_1, c_2 \in \mathbb{K}$ and a linear transformation $Z: V \to V$ we have $[X, c_1 Y + c_2 Z] = c_1 [X,Y] + c_2 [X,Z]$ and $[c_1 X + c_2 Y, Z] = c_1 [X,Z] + c_2 [Y,Z]$.    
\end{note} 

We present the following lemma in order to motivate the definitions of bidiagonal pair, triple, and triad.

\begin{lemma}
\cite[Lemma 2.2]{Funk-Neubauer17}
\label{thm:raise}
With reference to Notation \ref{note1} the following (i),(ii) hold.
\begin{enumerate}
\item[\rm (i)] Suppose that there exists an ordering $\lbrace Y_i \rbrace _{i=0}^d$ of the eigenspaces of $Y$ with \\ $X Y_i \subseteq Y_i + Y_{i+1} \, (0 \leq i \leq d)$, where $Y_{d+1} = 0$.  Then for $0 \leq i \leq d$, the restriction $[X, Y] |_{Y_i}$  maps $Y_i$ into $Y_{i+1}$.
\item[\rm (ii)] Suppose that there exists an ordering $\lbrace Y_i \rbrace _{i=0}^d$ of the eigenspaces of $Y$ with \\ $X Y_i \subseteq Y_{i-1} + Y_i \, (0 \leq i \leq d)$, where $Y_{-1} = 0$.  Then for $0 \leq i \leq d$, the restriction $[X, Y] |_{Y_i}$  maps $Y_i$ into $Y_{i-1}$.
\end{enumerate}
\end{lemma}

\begin{definition}
\rm
\cite[Definition 2.2]{Funk-Neubauer13}
\label{def:bdpair}
A {\it bidiagonal pair (BD pair) on $V$} is an ordered pair of linear transformations $A:V \to V$ and $A' :V \to V$ that satisfy the following three conditions.
\begin{enumerate}
\item[\rm (i)] Each of $A,\, A'$ is diagonalizable.
\item[\rm (ii)] There exists an ordering $\lbrace V_i \rbrace _{i=0}^d$ (resp.~$\lbrace V'_i \rbrace _{i=0}^D$) of the eigenspaces of $A$ (resp.~$A'$) with
\begin{align}
\label{eq:bp1}
A' V_i \subseteq V_i + V_{i+1} \qquad (0 \leq i \leq d), \\
\label{eq:bp2}
A V'_i \subseteq V'_i + V'_{i+1} \qquad (0 \leq i \leq D),
\end{align}
where each of $V_{d+1} , \, V'_{D+1}$ is equal to $0$.
\item[\rm (iii)] The restrictions
\begin{align}
\label{eq:bp3}
[A', A]^{d-2i} |_{V_i} : V_i \rightarrow V_{d-i} \qquad (0 \leq i \leq d/2), \\
\label{eq:bp4}
[A, A']^{D-2i} |_{V'_i} : V'_i \rightarrow V'_{D-i} \qquad (0 \leq i \leq D/2),
\end{align}
are bijections.
\end{enumerate}
\end{definition}

\begin{lemma}
\cite[Lemma 2.4]{Funk-Neubauer13}
\label{thm:dequalsdelta}
With reference to Definition \ref{def:bdpair}, we have $d=D$.
\end{lemma}

\begin{definition}
\rm
\cite[Definition 2.5]{Funk-Neubauer17}
\label{def:bdtriple}
A {\it bidiagonal triple (BD triple) on $V$} is an ordered triple of linear transformations $A:V\to V$, $A' :V \to V$, and $A'' :V \to V$ that satisfy the following three conditions.
\begin{enumerate}
\item[\rm (i)] Each of $A,\, A', \, A''$ is diagonalizable.
\item[\rm (ii)] There exists an ordering $\lbrace V_i \rbrace _{i=0}^d$ (resp.~$\lbrace V'_i \rbrace _{i=0}^{D}$) (resp.~$\lbrace V''_i \rbrace _{i=0}^{\delta}$) of the eigenspaces of $A$ (resp.~$A'$) (resp. ~$A''$) with
\begin{align}
\label{eq:bt1}
A' V_i \subseteq V_i + V_{i+1},  \qquad A'' V_i \subseteq V_{i-1} + V_i \qquad (0 \leq i \leq d), \\
\label{eq:bt2}
A'' V'_i \subseteq V'_i + V'_{i+1},  \qquad A V'_i \subseteq V'_{i-1} + V'_i \qquad (0 \leq i \leq D), \\
\label{eq:bt3}
AV''_i \subseteq V''_i + V''_{i+1},  \qquad A' V''_i \subseteq V''_{i-1} + V''_i \qquad (0 \leq i \leq \delta),
\end{align}
where each of $V_{d+1}, \, V_{-1},\, V'_{D+1}, \, V'_{-1}, \, V''_{\delta + 1}, \, V''_{-1}$ is equal to $0$.
\item[\rm (iii)] The restrictions
\begin{align}
\label{eq:bt4}
[A', A]^{d-2i} |_{V_i} : V_i \rightarrow V_{d-i}, \qquad  [A'', A]^{d-2i} |_{V_{d-i}} : V_{d-i} \rightarrow V_i \qquad (0 \leq i \leq d/2), \\
\label{eq:bt5}
[A'', A']^{D-2i} |_{V'_i} : V'_i \rightarrow V'_{D-i}, \qquad  [A, A']^{D-2i} |_{V'_{D-i}} : V'_{D-i} \rightarrow V'_i \,\,\,\,\,\, (0 \leq i \leq D/2), \\
\label{eq:bt6}
[A, A'']^{\delta-2i} |_{V''_i} : V''_i \rightarrow V''_{\delta-i}, \qquad  [A', A'']^{\delta-2i} |_{V''_{\delta-i}} : V''_{\delta-i} \rightarrow V''_i  \qquad (0 \leq i \leq \delta/2),
\end{align}
are bijections.
\end{enumerate}
\end{definition}

\begin{lemma}
\cite[Lemma 2.6]{Funk-Neubauer17}
\label{thm:dDdelta}
With reference to Definition \ref{def:bdtriple}, we have $d=D=\delta$.
\end{lemma}

\begin{definition}
\rm
\label{def:bdtriad}
A {\it bidiagonal triad (BD triad) on $V$} is an ordered triple of linear transformations $A:V\to V$, $A' :V \to V$, and $A'' :V \to V$ that satisfy the following three conditions.
\begin{enumerate}
\item[\rm (i)] Each of $A,\, A', \, A''$ is diagonalizable.
\item[\rm (ii)] There exists an ordering $\lbrace V_i \rbrace _{i=0}^d$ (resp.~$\lbrace V'_i \rbrace _{i=0}^{D}$) (resp.~$\lbrace V''_i \rbrace _{i=0}^{\delta}$) of the eigenspaces of $A$ (resp.~$A'$) (resp. ~$A''$) with
\begin{align}
\label{eq:btd1}
A' V_i \subseteq V_i + V_{i+1},  \qquad A'' V_i \subseteq V_i + V_{i+1} \qquad (0 \leq i \leq d), \\
\label{eq:btd2}
A'' V'_i \subseteq V'_i + V'_{i+1},  \qquad A V'_i \subseteq V'_i + V'_{i+1} \qquad (0 \leq i \leq D), \\
\label{eq:btd3}
AV''_i \subseteq V''_i + V''_{i+1},  \qquad A' V''_i \subseteq V''_i + V''_{i+1} \qquad (0 \leq i \leq \delta),
\end{align}
where each of $V_{d+1}, \, V'_{D+1}, \, V''_{\delta + 1}$ is equal to $0$.
\item[\rm (iii)] The restrictions
\begin{align}
\label{eq:btd4}
[A', A]^{d-2i} |_{V_i} : V_i \rightarrow V_{d-i}, \qquad  [A'', A]^{d-2i} |_{V_i} : V_i \rightarrow V_{d-i} \qquad (0 \leq i \leq d/2), \\
\label{eq:btd5}
[A'', A']^{D-2i} |_{V'_i} : V'_i \rightarrow V'_{D-i}, \qquad  [A, A']^{D-2i} |_{V'_i} : V'_i \rightarrow V'_{D-i} \,\,\,\,\,\, (0 \leq i \leq D/2), \\
\label{eq:btd6}
[A, A'']^{\delta-2i} |_{V''_i} : V''_i \rightarrow V''_{\delta-i}, \qquad  [A', A'']^{\delta-2i} |_{V''_i} : V''_i \rightarrow V''_{\delta-i}  \qquad (0 \leq i \leq \delta/2),
\end{align}
are bijections.
\end{enumerate}
\end{definition}

\begin{lemma}
\label{thm:dDdelta2}
With reference to Definition \ref{def:bdtriad}, we have $d=D=\delta$.
\end{lemma}

\noindent
{\it Proof:}  Let $A, \, A', \, A''$ denote a BD triad.  Then $A, \, A'$ is a BD pair.  So by Lemma \ref{thm:dequalsdelta} we have $d=D$.  Also, $A', \, A''$ is a BD pair.  So by Lemma \ref{thm:dequalsdelta} we have $D=\delta$.
\hfill $\Box$ \\

With reference to Definition \ref{def:bdtriad} we call $V$ the {\it vector space underlying} $A, \, A', \, A''$.  We say $A, \, A', \, A''$ is {\it over} $\K$.  In view of Lemma \ref{thm:dDdelta2}, for the remainder of this paper we use $d$ to index the eigenspaces of $A$, $A'$, and $A''$.  We call $d$ the {\it diameter} of $A, \, A', \, A''$.  The {\it diameter} of a BD pair and BD triple are defined similarly.    \\

For the remainder of this paper we assume that the field $\K$ is algebraically closed and characteristic zero. \\

Let $A, \, A', \, A''$ denote a BD triad of diameter $d$.  An ordering of the eigenspaces of $A$ (resp.~$A'$) (resp.~$A''$) is called {\it standard} whenever this ordering satisfies (\ref{eq:btd1}) (resp.~(\ref{eq:btd2})) (resp.~(\ref{eq:btd3})).  Let $\lbrace V_i \rbrace _{i=0}^d$ (resp.~$\lbrace V'_i \rbrace _{i=0}^d$) (resp.~$\lbrace V''_i \rbrace _{i=0}^d$) denote a standard ordering of the eigenspaces of $A$ (resp.~$A'$) (resp.~$A''$).  Then no other ordering of these eigenspaces is standard, and so the BD triad $A, \, A', \, A''$ uniquely determines these three standard orderings.  For $0 \leq i \leq d$, let $\theta_i$ (resp.~$\theta'_i $) (resp.~$\theta''_i $) denote the eigenvalue of $A$ (resp.~$A'$) (resp.~$A''$) corresponding to $V_i$ (resp.~$V'_i$) (resp.~$V''_i$).  We call $\lbrace \theta_i \rbrace _{i=0}^d$ (resp.~$\lbrace \theta'_i \rbrace _{i=0}^d$) (resp.~$\lbrace \theta''_i \rbrace _{i=0}^d$) the {\it first} (resp.~{\it second}) (resp.~{\it third}) {\it eigenvalue sequence} of $A, \, A', \, A''$.  For the remainder of this section, we adopt the notation from this paragraph.  All the terms defined in this paragraph have analogs for BD pairs and BD triples, see \cite[Section 2]{Funk-Neubauer13} and \cite[Section 2]{Funk-Neubauer17}.

\begin{lemma}
\label{thm:samedim}
For $0 \leq i \leq d$, the spaces $V_i$, $V_{d-i}$, $V'_i$, $V'_{d-i}$, $V''_i$, $V''_{d-i}$ all have the same dimension.
\end{lemma}

\noindent
{\it Proof:}  Since $A, \, A'$ is a BD pair then, by \cite[Lemma 2.7]{Funk-Neubauer13}, the spaces $V_i$, $V_{d-i}$, $V'_i$, $V'_{d-i}$ all have the same dimension.  Since $A', \, A''$ is a BD pair then, by \cite[Lemma 2.7]{Funk-Neubauer13}, the spaces $V'_i$, $V'_{d-i}$, $V''_i$, $V''_{d-i}$ all have the same dimension.
\hfill $\Box$ \\

\begin{definition}
\rm
\label{def:shape}
With reference to Lemma \ref{thm:samedim}, for $0 \leq i \leq d$, let $\rho_i$ denote the common dimension of  $V_i$, $V_{d-i}$, $V'_i$, $V'_{d-i}$, $V''_i$, $V''_{d-i}$.  We refer to the sequence $\lbrace \rho_i \rbrace _{i=0}^d$ as the {\it shape} of $A, \, A', \, A''$.
\end{definition}

The {\it shape} of a BD pair and BD triple are defined similarly.  See \cite[Lemma 2.8]{Funk-Neubauer13} and \cite[Lemma 2.9]{Funk-Neubauer17}. \\

The following definitions will be used in stating two of our main results (see Theorem \ref{thm:BDtriadtet} and Theorem \ref{thm:reduce}).

\begin{definition}
\rm
\label{def:thin}
Let $\lbrace \rho_i \rbrace _{i=0}^d$ denote the shape of $A, \, A', \, A''$.  We say that $A, \, A', \, A''$ is {\it thin} whenever $\rho_i = 1$ for $0 \leq i \leq d$.
\end{definition}

Similarly, a BD pair or BD triple is called {\it thin} whenever its shape is $(1,1,1, \ldots, 1)$.

\begin{definition}
\rm
\label{def:ae}
Let $\lbrace \sigma_i \rbrace _{i=0}^d$ and $\lbrace \tau_i \rbrace _{i=0}^d$ each denote a sequence of scalars taken from $\K$.  Let $X:V\to V$ and $Y:V \to V$ denote linear transformations.  Let $B, \, B', \, B''$ denote a BD triad on $V$.
\begin{enumerate}
\item[\rm (i)] We say $\lbrace \sigma_i \rbrace _{i=0}^d$ is {\it affine equivalent} to $\lbrace \tau_i \rbrace _{i=0}^d$, denoted $\lbrace \sigma_i \rbrace _{i=0}^d \sim \lbrace \tau_i \rbrace _{i=0}^d$, whenever there exist $r, s$ in $\K$ with $r \neq 0$ such that $\sigma_i = r \tau_i + s$ ($0 \leq i \leq d$).
\item[\rm (ii)] We say $X$ is {\it affine equivalent} to $Y$, denoted $X \sim Y$, whenever there exist $r, s$ in $\K$ with $r \neq 0$ such that $X = r Y + sI$.
\item[\rm (iii)] We say $A, \, A', \, A''$ is {\it affine equivalent} to $B, \, B', \, B''$ whenever $A \sim B$, $A' \sim B'$, and $A'' \sim B''$.
\end{enumerate}
\end{definition}

\section{The Tetrahedron Algebra}

\begin{definition}
\rm
\cite[Definition 1.1]{HarTer07}
\label{def:tet}
Let $\boxtimes$ denote the Lie algebra over $\mathbb{K}$ that has generators $$\{ X_{ij} \, | \, i, j \in \mathbb{I}, i \neq j \} \qquad \mathbb{I} = \{0, 1, 2, 3 \}$$ and the following relations:
\begin{enumerate}
\item For distinct $i, j \in \mathbb{I}$, $$X_{ij} + X_{ji} = 0.$$
\item For mutually distinct $h, i, j \in \mathbb{I}$, $$[X_{hi}, X_{ij}] = 2X_{hi} + 2X_{ij}.$$
\item For mutually distinct $h, i, j, k \in \mathbb{I}$, $$[X_{hi}, [X_{hi}, [X_{hi}, X_{jk}]]] = 4[X_{hi}, X_{jk}].$$
\end{enumerate}
\end{definition}
We call $\boxtimes$ the {\it tetrahedron algebra} and pronounce the symbol $\boxtimes$ as "tet". \\

The relations in Definition \ref{def:tet}(iii) are the well-known {\it Dolan-Grady} relations.  See \cite{AhnShig, DateRoan, Davies90, Davies91, DolanGrady, KlishPly}.   \\

The following definition will be used in stating two of our main results (see Theorem \ref{thm:tetBDtriad} and Theorem \ref{thm:BDtriadtet}).

\begin{definition}
\rm
\label{def:cornertri}
We call $X_{ru}, X_{su}, X_{tu}$ a {\it corner triad in $\boxtimes$} whenever $X_{ru}, X_{su}, X_{tu}$ are mutually distinct generators of $\boxtimes$.
\end{definition}

The name "corner triad" comes from the fact that the generators $X_{ru}, X_{su}, X_{tu}$ correspond to the three edges of the tetrahedron that meet at the vertex labeled by $u$. \\

The following information concerning $\boxtimes$-modules will be useful in proving one of our main results (Theorem \ref{thm:tetBDtriad}). 

\begin{lemma}
\cite[Theorem 3.8]{Hartwig07}
\label{tetdiag}
Let $V$ denote a finite-dimensional irreducible $\boxtimes$-module.  Then the following (i), (ii) hold.
\begin{enumerate}
\item For distinct $r, s \in \mathbb{I}$ the generator $X_{rs}$ is diagonalizable on $V$.  
\item There exists an integer $d \geq 0$ such that for all distinct $r,s \in \mathbb{I}$ the eigenvalues for $X_{rs}$ on $V$ are $\{ d-2i\}_{i=0}^d$.
\end{enumerate}
\end{lemma}

\begin{definition}
\rm
\cite[Definition 3.9]{Hartwig07}
\label{tetdiameter}
Let $V$ denote a finite-dimensional irreducible $\boxtimes$-module.  By the {\it diameter} of $V$ we mean the nonnegative integer $d$ from Lemma \ref{tetdiag}(ii).
\end{definition}

\begin{notation}
\rm
\label{def:tetespaces}
Let $V$ denote a finite-dimensional irreducible $\boxtimes$-module of diameter $d$.  For distinct $r, s \in \mathbb{I}$ and $0 \leq i \leq d$ we let $V_{rs}(d-2i)$ denote the eigenspace of $X_{rs}$ on $V$ corresponding to eigenvalue $d-2i$.  That is, 
\begin{eqnarray*}
V_{rs}(d-2i) := \{ v \in V \, | \, X_{rs} v = (d-2i) v\}.
\end{eqnarray*}  
\end{notation}

Note that by Definition \ref{def:tet}(i) we have 
\begin{eqnarray}
\label{tetgenneg}
X_{sr} = - X_{rs},
\end{eqnarray}
and so 
\begin{eqnarray}
\label{tetespaces}
V_{rs}(d-2i) = V_{sr}(2i-d). 
\end{eqnarray}

We now discuss the Lie algebra $\mathfrak{sl}_2$ and its relationship to $\boxtimes$.

\begin{definition}
\rm
\label{def:usl2}
Let $\mathfrak{sl}_2$ denote the Lie algebra over $\K$ that has a basis $h, \, e, \, f$ and Lie bracket
\begin{eqnarray*}
[h,e] = 2e, \qquad
[h,f] = -2f, \qquad
[e,f] = h.
\end{eqnarray*}
\end{definition}

The following two lemmas give a description of all finite-dimensional $\mathfrak{sl}_2$-modules.

\begin{lemma}
\cite[Theorem 6.3]{Humph72}
\label{thm:compred}
Each finite-dimensional $\mathfrak{sl}_2$-module $V$ is completely reducible; this means that $V$ is a direct sum of irreducible $\mathfrak{sl}_2$-modules.
\end{lemma}

The finite-dimensional irreducible $\mathfrak{sl}_2$-modules are described as follows.

\begin{lemma}
\cite[Theorem 7.2]{Humph72}
\label{thm:sl2mods}
There exists a family of finite-dimensional irreducible $\mathfrak{sl}_2$-modules
\begin{eqnarray*}
V(d), \qquad  d=0,1,2,\ldots
\end{eqnarray*}
with the following properties:  $V(d)$ has a basis $\lbrace v_i \rbrace _{i=0}^d$ such that $h.v_i = (d-2i) v_i$ for $0 \leq i \leq d$, $f.v_i = (i+1) v_{i+1}$ for $0 \leq i \leq d$, where $v_{d+1}=0$, and $e.v_i = (d-i+1) v_{i-1}$ for $0 \leq i \leq d$, where $v_{-1}=0$.  Moreover, every finite-dimensional irreducible $\mathfrak{sl}_2$-module is isomorphic to exactly one of the modules $V(d)$.
\end{lemma}

\begin{lemma}
\cite[Lemma 3.7]{Funk-Neubauer13}
\label{thm:2pieces}
Let $V$ denote an $\mathfrak{sl}_2$-module with finite positive dimension.  Define
\begin{eqnarray*}
V_{\rm{even}} &:=& \operatorname{span}\{\,v \in V\, |\, h.v = i \, v, \,\, i \in \mathbb{Z}, \, i \, even \,\}, \\
V_{\rm{odd}} &:=& \operatorname{span}\{\,v \in V\, |\, h.v = i \, v, \,\, i \in \mathbb{Z}, \, i \, odd \,\}.
\end{eqnarray*}
Then $V_{\rm{even}}$ and $V_{\rm{odd}}$ are $\mathfrak{sl}_2$-modules, and $V = V_{\rm{even}} + V_{\rm{odd}}$ (direct sum).
\end{lemma}

\begin{definition}
\rm
\cite[Definition 3.8]{Funk-Neubauer13}
\label{def:seg}
Let $V$ denote an $\mathfrak{sl}_2$-module with finite positive dimension.  With reference to Lemma \ref{thm:2pieces}, we say $V$ is {\it segregated} whenever $V=V_{\rm{even}}$ or $V=V_{\rm{odd}}$.
\end{definition}

The following two lemmas explain the relationship between $\mathfrak{sl}_2$ and $\boxtimes$.

\begin{lemma}
\cite[Lemma 3.2]{HarTer07}
\label{thm:ueq}
The Lie algebra $\mathfrak{sl}_2$ is isomorphic to the Lie algebra over $\K$ that has basis $X,\, Y, \, Z$ and Lie bracket
\begin{eqnarray*}
[X,Y] = 2X + 2Y, \qquad
[Y,Z] = 2Y + 2Z, \qquad
[Z,X] =2Z + 2X.
\end{eqnarray*}
An isomorphism with the presentation in Definition \ref{def:usl2} is given by:
\begin{eqnarray*}
X \rightarrow 2e - h, \qquad Y \rightarrow -2f - h, \qquad Z \rightarrow h.
\end{eqnarray*}
The inverse of this isomorphism is given by:
\begin{eqnarray*}
e \rightarrow (X + Z)/2, \qquad f \rightarrow -(Y + Z)/2, \qquad h \rightarrow Z.
\end{eqnarray*}
\end{lemma}

\begin{definition}
\rm
\cite[Definition 3.3]{Funk-Neubauer17}
\label{def:equitpair}
Let $X, \, Y, \, Z$ denote a basis for $\mathfrak{sl}_2$ satisfying the relations from Lemma \ref{thm:ueq}.  We refer to the ordered triple $X, \, Y, \, Z$ as an {\it equitable triple in $\mathfrak{sl}_2$}. 
\end{definition}

\begin{lemma}
\cite[Proposition 3.6]{HarTer07}
\label{sl2tetconnect}
Let $h, i, j$ denote mutually distinct elements of $\mathbb{I}$.  Then there exists a unique Lie algebra homomorphism from $\mathfrak{sl}_2$ to $\boxtimes$ that sends
\begin{eqnarray*}
X \rightarrow X_{hi}, \qquad Y \rightarrow X_{ij}, \qquad Z \rightarrow X_{jh}.
\end{eqnarray*}
\end{lemma}

\section{The Main Theorems}

The four theorems in this section make up the main conclusions of the paper.

\begin{theorem}
\label{thm:recur}
Let $A, \, A', \, A''$ denote a BD triad on $V$ of diameter $d$.  Let $\lbrace \theta_i \rbrace _{i=0}^d$, $\lbrace \theta'_i \rbrace _{i=0}^d$, and $\lbrace \theta''_i \rbrace _{i=0}^d$ denote the first, second, and third eigenvalue sequences of $A, \, A', \, A''$ respectively.  Suppose that $d \geq 2$.  Then for $1 \leq i \leq d-1$ we have
\begin{eqnarray}
\label{recur}
\frac{\theta_{i+1} - \theta_i}{\theta_i - \theta_{i-1}} = \frac{\theta'_{i+1} - \theta'_i}{\theta'_i - \theta'_{i-1}} = \frac{\theta''_{i+1} - \theta''_i}{\theta''_i - \theta''_{i-1}} = 1.
\end{eqnarray}
\end{theorem}

We refer to Theorem \ref{thm:recur} as the recurrence theorem.  The following definition will be used to state the next three theorems.

\begin{definition}
\label{def:reduced}
\rm
Let $A, \, A', \, A''$ denote a BD triad of diameter $d$.  We say $A, \, A', \, A''$ is {\it reduced} if the first, second, and third eigenvalue sequences of $A, \, A', \, A''$ are each $\lbrace 2i-d \rbrace _{i=0}^d$.
\end{definition}

See \cite[Definition 5.7]{Funk-Neubauer13} (resp. \cite[Definition 4.5]{Funk-Neubauer17}) for the corresponding definition of a reduced BD pair (resp. BD triple).  

\begin{theorem}
\label{thm:reduce}
Every BD triad is affine equivalent to a reduced BD triad.
\end{theorem}

We refer to Theorem \ref{thm:reduce} as the reducibility theorem.

\begin{note}
\label{reducedcases}
\rm
We make the following observation in order to motivate the next two theorems.  Let $X_{ru}, X_{su}, X_{tu}$ denote a corner triad in $\boxtimes$.  By Definition \ref{def:tet}(i),(ii) we have
\begin{eqnarray*}
[X_{ru}, X_{su}] &=& -2 X_{ru} + 2 X_{su}, \\
{[X_{su}, X_{tu}]} &=& -2 X_{su} + 2 X_{tu} \\
{[X_{tu}, X_{ru}]} &=& -2 X_{tu} + 2 X_{ru}.
\end{eqnarray*}
Let $A, \, A', \, A''$ denote a reduced BD triad.  It is shown in Lemma \ref{AA'A''brack} that
\begin{eqnarray*}
[A, A'] &=& -2 A + 2 A', \\
{[A', A'']} &=& -2 A' + 2 A'', \\
{[A'', A]} &=& -2 A'' + 2 A.
\end{eqnarray*}
\end{note}

The following two theorems provide a correspondence between reduced BD triads and irreducible $\boxtimes$-modules.

\begin{theorem}
\label{thm:tetBDtriad}
Let $V$ denote a finite-dimensional irreducible $\boxtimes$-module.  Then each corner triad in $\boxtimes$ acts on $V$ as a reduced BD triad.
\end{theorem}

We refer to the result in Theorem \ref{thm:tetBDtriad} as corner triads act as BD triads.  The next theorem can be thought of as a partial converse of Theorem \ref{thm:tetBDtriad}.

\begin{theorem}
\label{thm:BDtriadtet}
Let $A, \, A', \, A''$ denote a reduced thin BD triad on $V$.  Let $X_{ru}, X_{su}, X_{tu}$ denote a corner triad in $\boxtimes$.  Then there exists a $\boxtimes$-module structure on $V$ such that $(X_{ru} - A)V=0$, $(X_{su} - A')V=0$, and $(X_{tu} - A'')V=0$.  This $\boxtimes$-module structure on $V$ is irreducible.
\end{theorem}

We refer to the result in Theorem \ref{thm:BDtriadtet} as BD triads act as corner triads.  Theorem \ref{thm:reduce} and Theorem \ref{thm:BDtriadtet} together show that, in the thin case, every BD triad is affine equivalent to a BD triad of the type constructed in Theorem \ref{thm:tetBDtriad}.  \\

The $\boxtimes$-module constructed in Theorem \ref{thm:BDtriadtet} is an {\it evaluation module for $\boxtimes$} as defined in \cite[Definition 5.7]{ItoTerinpress4}.  See also \cite[Proposition 6.4]{ItoTerinpress4}.  Evaluation modules are the building blocks for all finite-dimensional irreducible $\boxtimes$-modules as described in \cite[Theorem 15.2]{ItoTerinpress4}. \\

In Theorem \ref{thm:BDtriadtet} the assumption that $A, \, A', \, A''$ is thin is needed to guarantee the conclusion of the theorem.  In Section 11 we provide a counterexample to demonstrate this.  

\section{Preliminaries}

In this section we develop some tools needed to prove our main results.

\begin{definition}
\label{def:brec}
\rm
Suppose that $d \geq 2$ and let $\lbrace \sigma_i \rbrace _{i=0}^d$ denote a sequence of distinct scalars taken from $\K$.  This sequence is called {\it $1$-recurrent} whenever
\begin{eqnarray*}
\frac{\sigma_{i+1} - \sigma_i}{\sigma_i - \sigma_{i-1}} = 1 \qquad (1 \leq i \leq d-1).
\end{eqnarray*}
\end{definition}

\begin{lemma}
\cite[Lemma 5.2]{Funk-Neubauer17}
\label{thm:brecuraffine}
Suppose that $d \geq 2$ and let $\lbrace \sigma_i \rbrace _{i=0}^d$ and $\lbrace \tau_i \rbrace _{i=0}^d$ each denote a sequence of distinct scalars taken from $\K$.  Assume $\lbrace \sigma_i \rbrace _{i=0}^d$ and $\lbrace \tau_i \rbrace _{i=0}^d$ are each $1$-recurrent.  Then $\lbrace \sigma_i \rbrace _{i=0}^d \sim \lbrace \tau_i \rbrace _{i=0}^d$.
\end{lemma}

\begin{lemma}
\cite[Lemma 5.3]{Funk-Neubauer17}
\label{thm:d<2}
Suppose that $d < 2$ and let $\lbrace \sigma_i \rbrace _{i=0}^d$ and $\lbrace \tau_i \rbrace _{i=0}^d$ each denote a sequence of distinct scalars taken from $\K$.  Then $\lbrace \sigma_i \rbrace _{i=0}^d \sim \lbrace \tau_i \rbrace _{i=0}^d$.
\end{lemma}

\begin{definition}
\rm
\cite[Definition 5.4]{Funk-Neubauer17}
\label{def:decomp}
A {\it decomposition} of $V$ is a sequence  $\lbrace U_i \rbrace _{i=0}^d$ consisting of nonzero subspaces of $V$ such that $V=\sum_{i=0}^{d} U_i$ (direct sum).  For any decomposition of $V$ we adopt the convention that $U_{-1} := 0$ and $U_{d+1} := 0$.
\end{definition}

\begin{lemma}
\cite[Lemma 5.5]{Funk-Neubauer17}
\label{thm:mapseq}
Let $\lbrace U_i \rbrace _{i=0}^d$ denote a decomposition of $V$.  Let $\lbrace \sigma_i \rbrace _{i=0}^d$ and $\lbrace \tau_i \rbrace _{i=0}^d$ each denote a sequence of distinct scalars taken from $\K$.  Let $X :V \to V$ (resp.~$Y: V \to V$) denote the linear transformation such that for $0 \leq i \leq d$, $U_i$ is an eigenspace for $X$ (resp.~$Y$) with eigenvalue $\sigma_i$ (resp.~$\tau_i$).  Then $\lbrace \sigma_i \rbrace _{i=0}^d \sim \lbrace \tau_i \rbrace _{i=0}^d$ if and only if $X \sim Y$.
\end{lemma}

We now develop some more properties of BD triads.  

\begin{lemma}
\label{thm:affineshift}
Let $A, \, A', \, A''$ denote a BD triad on $V$ of diameter $d$.  Let $\lbrace \theta_i \rbrace _{i=0}^d$, $\lbrace \theta'_i \rbrace _{i=0}^d$, $\lbrace \theta''_i \rbrace _{i=0}^d$ denote the first, second, third eigenvalue sequences of $A, \, A', \, A''$ respectively.  Let $r, \, s, \, t, \, u, \, v, \, w$ denote scalars in $\K$ with $r, \, t, \, v$ each nonzero.  Then $rA + sI , \, tA' + uI, \, vA'' + wI$ is a BD triad on $V$.  Moreover, $\lbrace r \theta_i + s \rbrace _{i=0}^d$, $\lbrace t \theta'_i  + u \rbrace _{i=0}^d$, $\lbrace v \theta''_i  + w \rbrace _{i=0}^d$ are the first, second, third eigenvalue sequences of $rA + sI , \, tA' + uI, \, vA'' + wI$ respectively.
\end{lemma}

\noindent
{\it Proof:}  Imitate the proof of  \cite[Lemma 2.10]{Funk-Neubauer13}.
\hfill $\Box$ \\

\begin{lemma}
\label{AA'A''brack}
Let $A, \, A', \, A''$ denote a reduced BD triad.  Then the following hold.
\begin{eqnarray}
\label{eq1}
[A, A'] &=& -2 A + 2 A', \\
\label{eq2}
{[A', A'']} &=& -2 A' + 2 A'', \\
\label{eq3}
{[A'', A]} &=& -2 A'' + 2 A.
\end{eqnarray}
\end{lemma}

\noindent
{\it Proof:}  We show (\ref{eq1}).  The proofs of (\ref{eq2}) and (\ref{eq3}) are similar.  Let $\lbrace \theta_i \rbrace _{i=0}^d$ (resp. $\lbrace \theta'_i \rbrace _{i=0}^d$) denote the first (resp. second) eigenvalue sequence of $A, \, A', \, A''$.  By Definition \ref{def:reduced} we have
\begin{eqnarray}
\label{eq4}
\theta_i = 2i-d = \theta'_i \qquad (0 \leq i \leq d).
\end{eqnarray}
First assume that $d = 0$.  Observe that $A = \theta_0 I$ and $A' = \theta'_0 I$.  By (\ref{eq4}), $\theta_0 =0$ and $\theta'_0 = 0$. Combining the previous two sentences we obtain (\ref{eq1}).  Now assume that $d \geq 1$.  Observe that $A, \, A'$ is a BD pair.  So by \cite[Theorem 5.3]{Funk-Neubauer13} there exist scalars $b, \alpha, \alpha^*, \gamma$ such that
\begin{eqnarray}
\label{eq5}
A A' - b A' A - \alpha A - \alpha^* A' - \gamma I = 0.
\end{eqnarray}
When $d=1$ then $b=1$ by \cite[Definition 5.5]{Funk-Neubauer13}.  When $d \geq 2$ then $b=1$ by (\ref{eq4}) and \cite[Lemma 9.1]{Funk-Neubauer13}.  Combining this with \cite[Theorem 8.1]{Funk-Neubauer13} and (\ref{eq4}) we have that $\alpha = -2, \alpha^* = 2, \gamma = 0$.  Substituting these values into (\ref{eq5}) we obtain (\ref{eq1}).
\hfill $\Box$ \\

\section{The Proofs of the Recurrence and Reducibility \\ Theorems}

In this section we prove Theorem \ref{thm:recur} and Theorem \ref{thm:reduce}.  First we prove the recurrence theorem. \\

\noindent
{\it Proof of Theorem \ref{thm:recur}:}   Let $A, \, A', \, A''$ denote a BD triad on $V$ of diameter $d$.  Let $\lbrace \theta_i \rbrace _{i=0}^d$, $\lbrace \theta'_i \rbrace _{i=0}^d$, and $\lbrace \theta''_i \rbrace _{i=0}^d$ denote the first, second, and third eigenvalue sequences of $A, \, A', \, A''$ respectively.  Suppose that $d \geq 2$.  For $1 \leq i \leq d-1$ abbreviate $\frac{\theta_{i+1} - \theta_i}{\theta_i - \theta_{i-1}}$, $\frac{\theta'_{i+1} - \theta'_i}{\theta'_i - \theta'_{i-1}}$, and $\frac{\theta''_{i+1} - \theta''_i}{\theta''_i - \theta''_{i-1}}$ as $T_i$, $T'_i$, and $T''_i$ respectively.  By construction $\lbrace \theta_i \rbrace _{i=0}^d$ is a list of distinct scalars.  So $T_i \neq 0$.  Similarly, $T'_i \neq 0$ and $T''_i \neq 0$.  Observe that $A, \, A'$ is a BD pair.  So by \cite[Theorem 5.1]{Funk-Neubauer13} we have
\begin{eqnarray}
\label{recur1}
T_i = T'^{-1}_i \qquad (1 \leq i \leq d-1).
\end{eqnarray}
Observe that $A', \, A''$ is a BD pair.  So by \cite[Theorem 5.1]{Funk-Neubauer13} we have
\begin{eqnarray}
\label{recur2}
T'_i = T''^{-1}_i \qquad (1 \leq i \leq d-1).
\end{eqnarray}
Observe that $A'', \, A$ is a BD pair.  So by \cite[Theorem 5.1]{Funk-Neubauer13} we have
\begin{eqnarray}
\label{recur3}
T''_i = T^{-1}_i \qquad (1 \leq i \leq d-1).
\end{eqnarray}
Combining (\ref{recur1}), (\ref{recur2}), (\ref{recur3}) we have that $T_i = T^{-1}_i$.  So $T_i$ must equal $1$ or $-1$.  But $\lbrace \theta_i \rbrace _{i=0}^d$ is a list of distinct scalars so $T_i \neq -1$.  Thus $T_i = 1$ for $1 \leq i \leq d-1$.  Similarly $T'_i = 1$ and $T''_i = 1$ for $1 \leq i \leq d-1$.
\hfill $\Box$ \\

We now prove the reducibility theorem.  The following proof is essentially the same as the proof of \cite[Theorem 4.6]{Funk-Neubauer17}.  For the sake of completeness and accessibility we reproduce the argument here in full.   \\

\noindent
{\it Proof of Theorem \ref{thm:reduce}:}  Let $A, \, A', \, A''$ denote a BD triad on $V$ of diameter $d$.  Let $\lbrace \theta_i \rbrace _{i=0}^d$, $\lbrace \theta'_i \rbrace _{i=0}^d$, and $\lbrace \theta''_i \rbrace _{i=0}^d$ denote the first, second, and third eigenvalue sequences of $A, \, A', \, A''$ respectively.  Let $\lbrace V_i \rbrace _{i=0}^d$, $\lbrace V'_i \rbrace _{i=0}^d$, and $\lbrace V''_i \rbrace _{i=0}^d$ denote the standard orderings of the eigenspaces of $A, \, A', \, A''$ respectively.  First we show that
\begin{eqnarray}
\label{reduce1}
{\rm \,\, each \,\, of \,\,} \lbrace \theta_i \rbrace _{i=0}^d, \lbrace \theta'_i \rbrace _{i=0}^d, \lbrace \theta''_i \rbrace _{i=0}^d {\rm \,\, is \,\, affine \,\, equivalent \,\, to \,\,} \lbrace 2i-d \rbrace _{i=0}^d.
\end{eqnarray}
If $d < 2$ then (\ref{reduce1}) is immediate from Lemma \ref{thm:d<2}.  Now assume that $d \geq 2$.  By Theorem \ref{thm:recur}, each of the sequences $\lbrace \theta_i \rbrace _{i=0}^d$, $\lbrace \theta'_i \rbrace _{i=0}^d$, $\lbrace \theta''_i \rbrace _{i=0}^d$ is $1$-recurrent.  Observe that the sequence $\lbrace 2i-d \rbrace _{i=0}^d$ is also $1$-recurrent.  So by Lemma \ref{thm:brecuraffine} we obtain (\ref{reduce1}).  Let $B :V \to V$ (resp.~$B' :V \to V$) (resp.~$B'' :V \to V$) be the linear transformation such that for $0 \leq i \leq d$, $V_i$ (resp.~$V'_i$) (resp.~$V''_i$) is an eigenspace for $B$ (resp.~$B'$) (resp.~$B''$) with eigenvalue $2i-d$.  From (\ref{reduce1}) and Lemma \ref{thm:mapseq} we have $A \sim B$, $A' \sim B'$, and $A'' \sim B''$.  Combining this with Lemma \ref{thm:affineshift} we find that $B, \, B', \, B''$ is a BD triad on $V$ of diameter $d$.  Also $A, \, A', \, A''$ is affine equivalent to $B, \, B', \, B''$.  Lemma \ref{thm:mapseq} and Lemma \ref{thm:affineshift} also show that each of the first, second, and third eigenvalue sequences of $B, \, B', \, B''$ is $\lbrace 2i-d \rbrace _{i=0}^d$.  Hence, $B, \, B', \, B''$ is reduced.  So $A, \, A', \, A''$ is affine equivalent to a reduced BD triad.
\hfill $\Box$ \\

\section{The Proof that Corner Triads Act as BD Triads}

In this section we prove Theorem \ref{thm:tetBDtriad}.  

\begin{note}
\label{basedefined}
\rm
Throughout the remainder of the paper we will refer to the {\it base} of a BD pair and BD triple.  See \cite[Definition 5.5]{Funk-Neubauer13} and \cite[Definition 2.13]{Funk-Neubauer17} for the definitions of base.   
\end{note}

\noindent
{\it Proof of Theorem \ref{thm:tetBDtriad}:}  Let $V$ denote a finite-dimensional irreducible $\boxtimes$-module of diameter $d$.  Let $X_{ru}, X_{su}, X_{tu}$ denote a corner triad in $\boxtimes$.  We show $X_{ru}, X_{su}, X_{tu}$ acts on $V$ as a reduced BD triad.  First, observe that by Lemma \ref{tetdiag} 
\begin{eqnarray}
\label{cortri1}
{\rm each \,\, of \,\,} X_{ru}, X_{su}, X_{tu} {\rm \,\, is \,\, diagonalizable \,\, on \,\,} V.
\end{eqnarray} 
By Lemma \ref{sl2tetconnect}, there exists an $\mathfrak{sl}_2$-module structure on $V$ such that $(X-X_{tu})V=0, (Y-X_{us})V=0, (Z-X_{st})V=0$ where $X, Y, Z$ is an equitable triple in $\mathfrak{sl}_2$.  We now check that this $\mathfrak{sl}_2$-module structure on $V$ is segregated.  By Lemma \ref{tetdiag} the action of $X_{st}$ on $V$ is diagonalizable with eigenvalues $\{ d-2i\}_{i=0}^d$.  Using the $\mathfrak{sl}_2$ isomorphism from Lemma \ref{thm:ueq} we have $(h-Z)V=0$.  Thus, $(h - X_{st})V = 0$ and so the action of $h$ on $V$ has eigenvalues $\{ d-2i\}_{i=0}^d$.  Hence, if $d$ is even (resp. odd) then $V = V_{\rm{even}}$ (resp. $V = V_{\rm{odd}}$).  This shows the $\mathfrak{sl}_2$-module structure on $V$ is segregated.  Applying \cite[Theorem 5.10]{Funk-Neubauer13} we have that $X_{tu}, X_{us}$ acts on $V$ as a reduced BD pair of diameter $d$ and base $1$.  By \cite[Definition 5.7]{Funk-Neubauer13} the eigenvalue (resp. dual eigenvalue) sequence of $X_{tu}, X_{us}$ is $\{ 2i-d \}_{i=0}^d$ (resp. $\{ d-2i\}_{i=0}^d$).  See \cite[Section 2]{Funk-Neubauer13} for the definitions of eigenvalue sequence and dual eigenvalue sequence of a BD pair.  From this, Definition \ref{def:bdpair}(ii), and Notation \ref{def:tetespaces} we have  
\begin{eqnarray*}
X_{tu} V_{us}(d-2i) \subseteq V_{us}(d-2i) + V_{us}(d-2(i+1)) \qquad (0 \leq i \leq d), \\
X_{us} V_{tu}(2i-d) \subseteq V_{tu}(2i-d) + V_{tu}(2(i+1)-d) \qquad (0 \leq i \leq d).
\end{eqnarray*}  
Combining the previous sentence with (\ref{tetgenneg}) and (\ref{tetespaces}) we have 
\begin{eqnarray}
\label{cortri2}
X_{tu} V_{su}(2i-d) \subseteq V_{su}(2i-d) + V_{su}(2(i+1)-d) \qquad (0 \leq i \leq d), \\
\label{cortri3}
X_{su} V_{tu}(2i-d) \subseteq V_{tu}(2i-d) + V_{tu}(2(i+1)-d) \qquad (0 \leq i \leq d).
\end{eqnarray} 
From Definition \ref{def:bdpair}(iii) and Notation \ref{def:tetespaces} we have 
\begin{eqnarray*}
{\rm the \,\, restriction \,\,} [X_{tu}, X_{us}]^{d-2i}: V_{us}(d-2i) \rightarrow V_{us}(2i-d) {\rm \,\, is \,\, a \,\, bijection \,\, for \,\,} 0 \leq i \leq d/2, \\
{\rm the \,\, restriction \,\,} [X_{tu}, X_{us}]^{d-2i}: V_{tu}(2i-d) \rightarrow V_{tu}(d-2i) {\rm \,\, is \,\, a \,\, bijection \,\, for \,\,} 0 \leq i \leq d/2. 
\end{eqnarray*}
Combining the previous sentence with (\ref{tetgenneg}) and (\ref{tetespaces}) we have 
\begin{eqnarray}
\label{cortri4}
{\rm the \,\, restriction \,\,} [X_{tu}, X_{su}]^{d-2i}: V_{su}(2i-d) \rightarrow V_{su}(d-2i) {\rm \,\, is \,\, a \,\, bijection \,\, for \,\,} 0 \leq i \leq d/2, \\
\label{cortri5}
{\rm the \,\, restriction \,\,} [X_{su}, X_{tu}]^{d-2i}: V_{tu}(2i-d) \rightarrow V_{tu}(d-2i) {\rm \,\, is \,\, a \,\, bijection \,\, for \,\,} 0 \leq i \leq d/2. 
\end{eqnarray}

By Lemma \ref{sl2tetconnect}, there exists an $\mathfrak{sl}_2$-module structure on $V$ such that $(X-X_{su})V=0, (Y-X_{ur})V=0, (Z-X_{rs})V=0$ where $X, Y, Z$ is an equitable triple in $\mathfrak{sl}_2$.  Repeating the argument from the previous paragraph in this case gives  
\begin{eqnarray}
\label{cortri6}
X_{su} V_{ru}(2i-d) \subseteq V_{ru}(2i-d) + V_{ru}(2(i+1)-d) \qquad (0 \leq i \leq d), \\
\label{cortri7}
X_{ru} V_{su}(2i-d) \subseteq V_{su}(2i-d) + V_{su}(2(i+1)-d) \qquad (0 \leq i \leq d), \\
\label{cortri8}
{\rm the \,\, restriction \,\,} [X_{su}, X_{ru}]^{d-2i}: V_{ru}(2i-d) \rightarrow V_{ru}(d-2i) {\rm \,\, is \,\, a \,\, bijection \,\, for \,\,} 0 \leq i \leq d/2, \\
\label{cortri9}
{\rm the \,\, restriction \,\,} [X_{ru}, X_{su}]^{d-2i}: V_{su}(2i-d) \rightarrow V_{su}(d-2i) {\rm \,\, is \,\, a \,\, bijection \,\, for \,\,} 0 \leq i \leq d/2. 
\end{eqnarray} 

By Lemma \ref{sl2tetconnect}, there exists an $\mathfrak{sl}_2$-module structure on $V$ such that $(X-X_{tu})V=0, (Y-X_{ur})V=0, (Z-X_{rt})V=0$ where $X, Y, Z$ is an equitable triple in $\mathfrak{sl}_2$.  Repeating the argument from the first paragraph again in this case gives  
\begin{eqnarray}
\label{cortri10}
X_{tu} V_{ru}(2i-d) \subseteq V_{ru}(2i-d) + V_{ru}(2(i+1)-d) \qquad (0 \leq i \leq d), \\
\label{cortri11}
X_{ru} V_{tu}(2i-d) \subseteq V_{tu}(2i-d) + V_{tu}(2(i+1)-d) \qquad (0 \leq i \leq d), \\
\label{cortri12}
{\rm the \,\, restriction \,\,} [X_{tu}, X_{ru}]^{d-2i}: V_{ru}(2i-d) \rightarrow V_{ru}(d-2i) {\rm \,\, is \,\, a \,\, bijection \,\, for \,\,} 0 \leq i \leq d/2, \\
\label{cortri13}
{\rm the \,\, restriction \,\,} [X_{ru}, X_{tu}]^{d-2i}: V_{tu}(2i-d) \rightarrow V_{tu}(d-2i) {\rm \,\, is \,\, a \,\, bijection \,\, for \,\,} 0 \leq i \leq d/2. 
\end{eqnarray} 
Combining (\ref{cortri1}) -- (\ref{cortri13}) we have that $X_{ru}, X_{su}, X_{tu}$ acts on $V$ as a BD triad.  It is immediate from (\ref{cortri2}), (\ref{cortri3}), (\ref{cortri6}),  (\ref{cortri7}), (\ref{cortri10}), (\ref{cortri11}) that the first, second, and third eigenvalue sequences of $X_{ru}, X_{su}, X_{tu}$ are each $\lbrace 2i-d \rbrace _{i=0}^d$.  Therefore, the BD triad $X_{ru}, X_{su}, X_{tu}$ is reduced.
\hfill $\Box$ \\

\section{The Raising Maps $R$ and $r$}

This section and the next two are devoted to the proof of Theorem \ref{thm:BDtriadtet}.  Throughout the remainder of the paper we will refer to the following assumption.

\begin{assumption}
\label{assump}
Let $A, \, A', \, A''$ denote a reduced BD triad on $V$ of diameter $d$.  Let $\lbrace V_i \rbrace _{i=0}^d$, $\lbrace V'_i \rbrace _{i=0}^d$, and $\lbrace V''_i \rbrace _{i=0}^d$ denote the standard orderings of the eigenspaces of $A, \, A', \, A''$ respectively. 
\end{assumption}

\begin{definition}
\label{defraise}
\rm
With reference to Assumption \ref{assump}, let $R: V \rightarrow V$ and $r: V \rightarrow V$ denote the following linear transformations: 
\begin{eqnarray*}
R := A - A'', \qquad \qquad r := A' - A''. 
\end{eqnarray*}
\end{definition}

\begin{note}
\label{Rrnotzero}
\rm
With reference to Assumption \ref{assump} note that by Definition \ref{def:bdtriad}(iii) $A \neq A''$, and so $R \neq 0$.  Similarly, $A' \neq A''$ and so $r \neq 0$.  Lastly, $A \neq A'$ and so $R \neq r$. 
\end{note}

We refer to the linear transformations from Definition \ref{defraise} as {\it raising maps}.  The following lemma explains this terminology.  

\begin{lemma}
\label{raise1}
With reference to Definition \ref{defraise} the following (i), (ii) hold.
\begin{enumerate}
\item[\rm (i)]  $R V''_i \subseteq V''_{i+1} \qquad (0 \leq i \leq d)$.
\item[\rm (ii)]  $r V''_i \subseteq V''_{i+1} \qquad (0 \leq i \leq d)$.
\end{enumerate}
\end{lemma}

\noindent
{\it Proof:} (i) Adopt Assumption \ref{assump}.  Observe that $A, \, -A''$ is a reduced BD pair of base $1$.  Combining this with \cite[Lemma 6.7(i)]{Funk-Neubauer13} we have, for $0 \leq i \leq d$ and $v \in V''_i$,
\begin{eqnarray*}
[A, -A'']v &=& 2 (A - (2i-d)I)v \\
&=& 2 (A - A'')v.
\end{eqnarray*}
Combining this with Lemma \ref{thm:raise}(i) and Definition \ref{defraise} we obtain (i).  \\
(ii) Similar to (i).
\hfill $\Box$ \\

\begin{lemma}
\label{raise2}
With reference to Definition \ref{defraise} we have $Rr = rR$.
\end{lemma}

\noindent
{\it Proof:}  Adopt Assumption \ref{assump}.  Observe that $A, \, -A'$ is a reduced BD pair of base $1$.  Combining this with \cite[Lemma 12.2]{Funk-Neubauer13} we have 
\begin{eqnarray}
\label{raise2a}
[A, -A'] = 2A - 2A'.
\end{eqnarray}
From Definition \ref{defraise} we have $A = R + A''$ and $A' = r + A''$.  Substituting these into (\ref{raise2a}) and simplifying gives 
\begin{eqnarray}
\label{raise2b}
[R, r] + [R, A''] + [A'', r] + 2R - 2r = 0.
\end{eqnarray}
Recall for $0 \leq i \leq d$ that $V''_i$ is the eigenspace for $A''$ corresponding to the eigenvalue $2i-d$.  Combining this with Lemma \ref{raise1} and (\ref{raise2b}) we have 
\begin{eqnarray*}
([R,r] + (2i-d)R - (2i+2-d)R + (2i+2-d)r - (2i-d)r + 2R - 2r)V''_i = 0.
\end{eqnarray*}  
Simplifying this gives $(Rr - rR)V''_i = 0$ for $0 \leq i \leq d$.  Since $A''$ is diagonalizable then $\lbrace V''_i \rbrace _{i=0}^d$ is a decomposition of $V$.  Therefore, $Rr - rR = 0$ and the result follows.  
\hfill $\Box$ \\

\begin{lemma}
\label{raise3}
Adopt Assumption \ref{assump} and assume that the BD triad $A, \, A', \, A''$ is thin.  With reference to Definition \ref{defraise} the following (i),(ii) hold.
\begin{enumerate}
\item[\rm (i)] The restriction $R|_{V''_i}: V''_i \rightarrow V''_{i+1}$ is a bijection for $0 \leq i \leq d-1$.
\item[\rm (ii)] The restriction $r|_{V''_i}: V''_i \rightarrow V''_{i+1}$ is a bijection for $0 \leq i \leq d-1$.
\end{enumerate}  
\end{lemma}

\noindent
{\it Proof:} (i) Observe that $A, \, -A''$ is a reduced BD pair of base $1$.  Combining this with \cite[Lemma 6.7(i)]{Funk-Neubauer13} and Definition \ref{defraise} we have
\begin{eqnarray}
\label{raise3a}
[A, -A'']|_{V''_i} = 2 R|_{V''_i} \qquad (0 \leq i \leq d).
\end{eqnarray}
Case 1:  $0 \leq i < d/2$.  Combining \cite[Lemma 2.5]{Funk-Neubauer13} with (\ref{raise3a}) shows that the restriction $R|_{V''_i}: V''_i \rightarrow V''_{i+1}$ is an injection.  Since  $A, \, A', \, A''$ is thin then ${\rm dim}(V''_i) = 1 = {\rm dim}(V''_{i+1})$.  Thus, the restriction $R|_{V''_i}: V''_i \rightarrow V''_{i+1}$ is also a surjection and hence a bijection. \\
Case 2:  $d/2 \leq i \leq d-1$.  Combining \cite[Lemma 2.5]{Funk-Neubauer13} with (\ref{raise3a}) shows that the restriction $R|_{V''_i}: V''_i \rightarrow V''_{i+1}$ is a surjection.  Since  $A, \, A', \, A''$ is thin then ${\rm dim}(V''_i) = 1 = {\rm dim}(V''_{i+1})$.  Thus, the restriction $R|_{V''_i}: V''_i \rightarrow V''_{i+1}$ is also an injection and hence a bijection. \\
(ii) Similar to (i).
\hfill $\Box$ \\

\begin{lemma}
\label{raise4}
Adopt Assumption \ref{assump} and assume that the BD triad $A, \, A', \, A''$ is thin.  With reference to Definition \ref{defraise} there exists a nonzero $c \in \mathbb{K}$ such that $r = cR$. 
\end{lemma}

\noindent
{\it Proof:}  Since $A''$ is diagonalizable then $\lbrace V''_i \rbrace _{i=0}^d$ is a decomposition of $V$.  So to prove the desired result it suffices to show 
\begin{eqnarray}
\label{raise4a}
{\rm there \,\, exists \,\, a \,\, nonzero \,\,} c \in \mathbb{K} {\rm \,\, such \,\, that \,\,} (r - cR)V''_i = 0 {\rm \,\, for \,\,} 0 \leq i \leq d. 
\end{eqnarray}
Recall that since $A, \, A', \, A''$ is thin we have
\begin{eqnarray}
\label{raise4b}
{\rm dim}(V''_i) = 1 \qquad 0 \leq i \leq d. 
\end{eqnarray}
We show (\ref{raise4a}) by induction on $i$.  First we show (\ref{raise4a}) holds for $i=0$.  Let $v$ denote a nonzero vector in $V''_0$.  Note that $\{ v \}$ is a basis for $V''_0$ by (\ref{raise4b}).  Observe that $Rv \in V''_1$ and $Rv \neq 0$ by Lemma \ref{raise3}(i).  Also $rv \in V''_1$ and $rv \neq 0$ by Lemma \ref{raise3}(ii).  Combining the previous two sentences with (\ref{raise4b}) shows there exists a nonzero $c \in \mathbb{K}$ such that $rv = cRv$.  This shows (\ref{raise4a}) holds for $i=0$ since $\{ v \}$ is a basis for $V''_0$.  Now let $i \geq 1$ and assume (\ref{raise4a}) holds for $i$.  Let $w$ denote a nonzero vector in $V''_i$.  Note that $Rw \in V''_{i+1}$ and $Rw \neq 0$ by Lemma \ref{raise3}(i).  Combining this with (\ref{raise4b}) shows $\{ Rw \}$ is a basis for $V''_{i+1}$.  By Lemma \ref{raise2} and the induction hypothesis we have $(r - cR) Rw = R(r - cR)w = 0$.  This shows  $(\ref{raise4a})$ holds for $i+1$ since $\{ Rw \}$ is a basis for $V''_{i+1}$. 
\hfill $\Box$ \\

\section{The Linear Transformations $B, \, B', \, B''$}

\begin{lemma}
\label{defB}
Adopt Assumption \ref{assump}.  Then there exists a linear transformation $B: V \rightarrow V$ such that $A', \, -A'', \, B$ is a reduced BD triple on $V$ of base $1$.  Moreover, if the BD triad $A, \, A', \, A''$ is thin then the BD triple $A', \, -A'', \, B$ is also thin.      
\end{lemma}

\noindent
{\it Proof:}  Observe that $A', \, -A''$ is a reduced BD pair on $V$ of base $1$.  So by \cite[Theorem 4.1]{Funk-Neubauer17} there exists a linear transformation $\overline{B}: V \rightarrow V$ such that $A', \, -A'', \, \overline{B}$ is a BD triple on $V$ of base $1$.  Moreover, $\{ 2i-d \}_{i=0}^d$ is both the first and second eigenvalue sequence of $A', \, -A'', \, \overline{B}$.  Let $\{ \overline{\beta_i} \}_{i=0}^d$ denote the third eigenvalue sequence of $A', \, -A'', \, \overline{B}$.  See \cite[Section 2]{Funk-Neubauer17} for the definition of first, second, third eigenvalue sequence of a BD triple.  By \cite[Lemma 2.12]{Funk-Neubauer17}, Lemma \ref{thm:brecuraffine}, and Lemma \ref{thm:d<2} we have $\lbrace 2i-d \rbrace _{i=0}^d \sim \lbrace \overline{\beta_i} \rbrace _{i=0}^d$.  So by Definition \ref{def:ae}(i) there exist $r, s$ in $\K$ with $r \neq 0$ such that $2i-d = r \overline{\beta_i} + s$ for $0 \leq i \leq d$.  Define the linear transformation $B: V \rightarrow V$ as $B := r \overline{B} + s I$.  Thus, by \cite[Lemma 5.7]{Funk-Neubauer17} we have that $A', \, -A'', \, B$ is a reduced BD triple on $V$ of base $1$.  Now suppose that the BD triad $A, \, A', \, A''$ is thin.  So both the BD triad $A, \, A', \, A''$ and the BD pair $A', \, -A''$ have shape $(1,1,1, \ldots, 1)$.  From this and \cite[Theorem 4.1]{Funk-Neubauer17} the BD triple $A', \, -A'', \, B$ also has shape $(1,1,1, \ldots, 1)$.  Thus, the BD triple $A', \, -A'', \, B$ is thin.        
\hfill $\Box$ \\

\begin{definition}
\label{defa}
\rm
Let $c$ denote the nonzero scalar in $\mathbb{K}$ from Lemma \ref{raise4}.  Define the scalar $a$ in $\mathbb{K}$ as $a := 1 - c^{-1}$.
\end{definition}
 
\begin{note}
\label{anot0not1}
\rm
Let $c$ denote the scalar from Lemma \ref{raise4}, and let $a$ denote the scalar from Definition \ref{defa}.  By Note \ref{Rrnotzero} and Lemma \ref{raise4} we have that $c \neq 1$.  Thus, $a \neq 0$.  Also, it is clear from Definition \ref{defa} that $a \neq 1$.  
\end{note} 
 
\begin{definition}
\label{defB'B''}
\rm
With reference to Lemma \ref{defB} and Definition \ref{defa}, let $B': V \rightarrow V$ and $B'': V \rightarrow V$ denote the following linear transformations: 
\begin{eqnarray}
\label{defB'}
B' &:=& (a^{-1} - 1)^{-1}A'' + (a-1)^{-1}B, \\ 
\label{defB''}
B'' &:=& (1-a^{-1})A' - a^{-1}B. 
\end{eqnarray}
\end{definition}
 
\begin{lemma}
\label{cornerrels}
Adopt Assumption \ref{assump} and assume that the BD triad $A, \, A', \, A''$ is thin.  With reference to Lemma \ref{defB} and Definition \ref{defB'B''} the following hold.
\begin{eqnarray}
\label{correl1}
A &=& (1-a)A' + aA'', \\
\label{correl2}
A' &=& (1-a^{-1})^{-1}A'' + (1-a)^{-1} A, \\
\label{correl3}
A'' &=& a^{-1}A + (1-a^{-1})A', \\
\label{correl4}
B &=& (a-1)B' + aA'', \\
\label{correl5}
A'' &=& a^{-1}B + (a^{-1} - 1)B', \\
\label{correl6}
B &=& -aB''+(a-1)A', \\
\label{correl7}
A' &=& (a-1)^{-1}B +(1-a^{-1})^{-1}B'', \\
\label{correl8}
A &=& (1-a)B' - aB'', \\
\label{correl9}
B' &=& (a^{-1}-1)^{-1}B'' + (1-a)^{-1}A, \\
\label{correl10}
B'' &=& -a^{-1}A + (a^{-1}-1)B'.
\end{eqnarray}  
\end{lemma}

\noindent
{\it Proof:}  By Definition \ref{defa} we have $c = (1-a)^{-1}$.  Combining this with Definition \ref{defraise} and Lemma \ref{raise4} gives $A' - A'' = (1-a)^{-1}(A - A'')$.  Solving this for  $A$ gives (\ref{correl1}).  Solving (\ref{correl1}) for $A'$ gives (\ref{correl2}).  Solving (\ref{correl1}) for $A''$ gives (\ref{correl3}).  Solving (\ref{defB'}) for $B$ gives (\ref{correl4}).  Solving (\ref{defB'}) for $A''$ gives (\ref{correl5}).  Solving (\ref{defB''}) for $B$ gives (\ref{correl6}).  Solving (\ref{defB''}) for $A'$ gives (\ref{correl7}).  Substituting (\ref{correl5}),(\ref{correl7}) into (\ref{correl1}) and simplifying gives (\ref{correl8}).  Solving (\ref{correl8}) for $B'$ gives (\ref{correl9}).  Solving (\ref{correl8}) for $B''$ gives (\ref{correl10}).            
\hfill $\Box$ \\

\begin{lemma}
\label{brack1}
Adopt Assumption \ref{assump}.  With reference to Lemma \ref{defB} the following hold.
\begin{eqnarray}
\label{brack1a}
[A'',B] &=& 2A'' - 2B, \\
\label{brack1b}
[B,A'] &=& 2B + 2A'.
\end{eqnarray}  
\end{lemma}

\noindent
{\it Proof:}  Immediate from Lemma \ref{defB} and \cite[Corollary 8.2(i)]{Funk-Neubauer17}. 
\hfill $\Box$ \\

\begin{lemma}
\label{brack2}
Adopt Assumption \ref{assump} and assume that the BD triad $A, \, A', \, A''$ is thin.  With reference to Lemma \ref{defB} and Definition \ref{defB'B''} the following hold.
\begin{eqnarray}
\label{brack2a}
[B',A''] &=& 2B' + 2A'', \\
\label{brack2b}
[B',B] &=& 2B' + 2B, \\
\label{brack2c}
[B,B''] &=& 2B + 2B'', \\
\label{brack2d}
[A',B''] &=& 2A' - 2B'', \\
\label{brack2e}
[B'',B'] &=& 2B'' + 2B', \\
\label{brack2f}
[B'',A] &=& 2B'' + 2A, \\
\label{brack2g}
[A,B'] &=& 2A - 2B'.
\end{eqnarray}  
\end{lemma}

\noindent
{\it Proof:}  Substituting (\ref{correl4}) into (\ref{brack1a}) and simplifying gives (\ref{brack2a}).  Substituting (\ref{correl5}) into (\ref{brack1a}) and simplifying gives (\ref{brack2b}).  Substituting (\ref{correl7}) into (\ref{brack1b}) and simplifying gives (\ref{brack2c}).  Substituting (\ref{correl6}) into (\ref{brack1b}) and simplifying gives (\ref{brack2d}).  Substituting (\ref{correl5}),(\ref{correl7}) into (\ref{eq2}) and simplifying using (\ref{brack2b}),(\ref{brack2c}) gives (\ref{brack2e}).  Substituting (\ref{correl3}),(\ref{correl6}) into (\ref{brack1a}) and simplifying using (\ref{eq1}),(\ref{brack2d}) gives (\ref{brack2f}).  Substituting (\ref{correl2}),(\ref{correl4}) into (\ref{brack1b}) and simplifying using (\ref{eq3}),(\ref{brack2a}) gives (\ref{brack2g}). 
\hfill $\Box$ \\

\begin{lemma}
\label{DG}
Adopt Assumption \ref{assump} and assume that the BD triad $A, \, A', \, A''$ is thin.  With reference to Lemma \ref{defB} and Definition \ref{defB'B''} the following (i)--(vi) hold.
\begin{enumerate}
\item[\rm (i)] $[A, [A, [A, B]]] = 4[A, B].$
\item[\rm (ii)] $[B, [B, [B, A]]] = 4[B, A].$
\item[\rm (iii)] $[A', [A', [A', B']]] = 4[A', B'].$
\item[\rm (iv)] $[B', [B', [B', A']]] = 4[B', A'].$
\item[\rm (v)] $[A'', [A'', [A'', B'']]] = 4[A'', B''].$
\item[\rm (vi)] $[B'', [B'', [B'', A'']]] = 4[B'', A''].$
\end{enumerate}  
\end{lemma}

\noindent
{\it Proof:} (i) We have 
\begin{eqnarray*}
[A, [A, [A, B]]] &=& [A, [A, [A, (a-1)B' + aA'']]] \qquad {\rm by \,\, (\ref{correl4})}  \\
&=& [A, [A, (a-1)(2A - 2B') + a(-2 A + 2 A'')]] \qquad {\rm by \,\, (\ref{eq3}),(\ref{brack2g})}  \\
&=& [A, - 2(a-1)[A,B'] + 2a[A, A'']] \\
&=& 4(a-1)[A,B'] + 4a[A,A'']  \qquad {\rm by \,\, (\ref{eq3}),(\ref{brack2g})} \\
&=& 4[A, B] \qquad {\rm by \,\, (\ref{correl4})}.
\end{eqnarray*} 
(ii) We have 
\begin{eqnarray*}
[B, [B, [B, A]]] &=& [B, [B, [B, (1-a)A' + aA'']]] \qquad {\rm by \,\, (\ref{correl1})}  \\
&=& [B, [B, (1-a)(2B + 2A') + a(2B - 2A'')]] \qquad {\rm by \,\, (\ref{brack1a}),(\ref{brack1b})}  \\
&=& [B, 2(1-a)[B,A'] - 2a[B, A'']] \\
&=& 4(1-a)[B,A'] + 4a[B,A'']  \qquad {\rm by \,\, (\ref{brack1a}),(\ref{brack1b})} \\
&=& 4[B, A] \qquad {\rm by \,\, (\ref{correl1})}.
\end{eqnarray*} 
(iii),(v) Similar to (i). \\
(iv),(vi) Similar to (ii).
\hfill $\Box$ \\

\section{The Proof that BD triads Act as Corner Triads}

In this section we prove Theorem \ref{thm:BDtriadtet}.  

\begin{lemma}
\label{RR1}
With reference to Notation \ref{note1} the following (i),(ii) are equivalent.
\begin{enumerate}
\item[\rm (i)] $[X,Y] = 2X + 2Y.$
\item[\rm (ii)] $[-Y,-X] = -2Y - 2X.$
\end{enumerate}
\end{lemma}

\noindent
{\it Proof:}  Immediate from Note \ref{brackprops}.
\hfill $\Box$ \\

\begin{lemma}
\label{RR2}
With reference to Notation \ref{note1} the following (i)--(iv) are equivalent.
\begin{enumerate}
\item[\rm (i)] $[X, [X, [X, Y]]] = 4[X, Y].$
\item[\rm (ii)] $[X, [X, [X, -Y]]] = 4[X, -Y].$
\item[\rm (iii)] $[-X, [-X, [-X, Y]]] = 4[-X, Y].$
\item[\rm (iv)] $[-X, [-X, [-X, -Y]]] = 4[-X, -Y].$
\end{enumerate}
\end{lemma}

\noindent
{\it Proof:}  Immediate from Note \ref{brackprops}.
\hfill $\Box$ \\

\begin{lemma}
\label{bijlem}
Let $X, \, Y, \, Z$ denote a thin BD triple on $V$.  Let $\lbrace Y_i \rbrace _{i=0}^d$ denote the standard ordering of the eigenspaces of $Y$.  Then the following (i),(ii) hold.    
\begin{enumerate}
\item[\rm (i)] The restriction $[X,Y]^j|_{Y_i}: Y_i \rightarrow Y_{i-j}$ is a bijection for $0 \leq i \leq d$ and $0 \leq j \leq i$.
\item[\rm (ii)] The restriction $[Z,Y]^j|_{Y_i}: Y_i \rightarrow Y_{i+j}$ is a bijection for $0 \leq i \leq d$ and $0 \leq j \leq d-i$.
\end{enumerate}
\end{lemma}

\noindent
{\it Proof:} (i) Recall by Lemma \ref{thm:raise}(ii) and (\ref{eq:bt2}) that the restriction $[X,Y]^j|_{Y_i}$ maps $Y_i$ to $Y_{i-j}$ for $0 \leq i \leq d$ and $0 \leq j \leq i$.  We now show that this mapping is a bijection.  Since $X, \, Y, \, Z$ is thin we have 
 \begin{eqnarray}
\label{bij1}
{\rm dim}(Y_i) = 1 {\rm \,\, for \,\,} 0 \leq i \leq d.
\end{eqnarray}
The result is trivially true for $i = 0$.  Now let $i \geq 1$.  Since the composition of two bijections is also a bijection it suffices to show that 
\begin{eqnarray}
\label{bij2}
{\rm The \,\, restriction \,\,} [X,Y]|_{Y_i}: Y_i \rightarrow Y_{i-1} {\rm \,\, is \,\, a \,\, bijection \,\, for \,\,} 1 \leq i \leq d.
\end{eqnarray}
Case 1:  $1 \leq i \leq d/2$.  First we show that the restriction $[X,Y]|_{Y_i}: Y_i \rightarrow Y_{i-1}$ is a surjection.  Let $v \in Y_{i-1}$.  By (\ref{eq:bt5}) there exists $w \in Y_{d-i+1}$ such that $[X,Y]^{d-2i+2}w = v$.  Observe that $[X,Y]^{d-2i+1}w \in Y_i$ and $[X,Y] [X,Y]^{d-2i+1}w = v$.  This shows that the restriction $[X,Y]|_{Y_i}: Y_i \rightarrow Y_{i-1}$ is a surjection.  From this and (\ref{bij1}) we have (\ref{bij2}). \\
Case 2:  $d/2 < i \leq d$.  First we show that the restriction $[X,Y]|_{Y_i}: Y_i \rightarrow Y_{i-1}$ is an injection.  Let $v \in {\rm ker}([X,Y]|_{Y_i})$.  So, $[X,Y]v = 0$ and thus $[X,Y]^{2i-d}v = 0$.  From this and (\ref{eq:bt5}) we have $v=0$.  This shows that the restriction $[X,Y]|_{Y_i}: Y_i \rightarrow Y_{i-1}$ is an injection.  From this and (\ref{bij1}) we have (\ref{bij2}).     \\
(ii) Similar to (i).
\hfill $\Box$ \\

\noindent
{\it Proof of Theorem \ref{thm:BDtriadtet}:}  Adopt Assumption \ref{assump} and assume that the BD triad $A, \, A', \, A''$ is thin.  Let $X_{ru}, X_{su}, X_{tu}$ denote a corner triad in $\boxtimes$.  Let $B$ denote the linear transformation from Lemma \ref{defB}.  Let $B', \, B''$ denote the linear transformations from Definition \ref{defB'B''}.  Comparing Definition \ref{def:tet} with Lemmas \ref{AA'A''brack}, \ref{brack1}, \ref{brack2}, \ref{DG}, \ref{RR1}, and \ref{RR2} we find there exists a $\boxtimes$-module structure on $V$ such that 
\begin{eqnarray*}
(X_{ru} - A)V=0, \qquad (X_{ur} + A)V=0, \qquad (X_{su} - A')V=0, \qquad (X_{us} + A')V=0, \\ 
(X_{tu} - A'')V=0, \qquad (X_{ut} + A'')V=0, \qquad (X_{ts} - B)V=0, \qquad (X_{st} + B)V=0, \\ 
(X_{rt} - B')V=0, \qquad (X_{tr} + B')V=0, \qquad (X_{sr} - B'')V=0, \qquad (X_{rs} + B'')V=0. 
\end{eqnarray*}
We now show that this $\boxtimes$-module structure on $V$ is irreducible.  Let $W$ denote a nonzero $\boxtimes$-submodule of $V$.  We show that $W=V$.  Let $0 \neq w \in W$.  Since $A''$ is diagonalizable then  $\lbrace V''_i \rbrace _{i=0}^d$ is a decomposition of $V$.  So for $0 \leq i \leq d$ there exists $v_i'' \in V_i''$ such that $w = \sum_{i=0}^d v_i''$.  Define $t := {\rm min}\{ i \, | \, 0 \leq i \leq d, \, v_i'' \neq 0\}$, and observe that 
\begin{eqnarray}
\label{sum1}
w = \sum_{i=t}^d v_i''.
\end{eqnarray}
Recall by Lemma \ref{defB} that $A', \, -A'', \, B$ is a thin BD triple and note that $\{ V_{d-i}''\}_{i=0}^d$ is the standard ordering of the eigenspaces of $-A''$.  By Lemma \ref{bijlem}(i) we have that the restriction $[A',-A'']^{d-t}|_{V_t''}: V_t'' \rightarrow V_d''$ is a bijection.  Abbreviate the vector $[A',-A'']^{d-t}v_t''$ as $v$.  Recall $0 \neq v_t'' \in V_t''$ by construction, and so 
\begin{eqnarray}
\label{vnotzero}
0 \neq v \in V_d''.
\end{eqnarray}  
By Lemma \ref{bijlem}(ii) we have that for $0 \leq r \leq d$ the restriction $[B,-A'']^r|_{V_d''}: V_d'' \rightarrow V_{d-r}''$ is a bijection.  Combining this with (\ref{vnotzero}) we have $0 \neq [B,-A'']^r v \in V_{d-r}''$ for $0 \leq r \leq d$.  Since the BD triple $A', \, -A'', \, B$ is thin then ${\rm dim}(V_{d-r}'') = 1$ for $0 \leq r \leq d$.  Combining the previous two sentences with the fact that $A''$ is diagonalizable shows 
\begin{eqnarray}
\label{basis}
\{ [B,-A'']^r v \}_{r=0}^d {\rm \,\, is \,\, a \,\, basis \,\, for \,\,} V.
\end{eqnarray}  
Since $W$ is a $\boxtimes$-submodule of $V$ then 
\begin{eqnarray}
\label{submod}
A'W \subseteq W, \qquad -A''W \subseteq W, \qquad BW \subseteq W.
\end{eqnarray}
By (\ref{sum1}) and Lemma \ref{thm:raise}(i) we have $[A',-A'']^{d-t}w = v$.  From this and (\ref{submod}) we have that $v \in W$.  Combining this with (\ref{basis}) and (\ref{submod}) shows that $W=V$.  Therefore, the $\boxtimes$-module structure on $V$ is irreducible.  
\hfill $\Box$ \\

\begin{note}
\label{BDtriplestriadspairs}
\rm
Adopt Assumption \ref{assump} and assume that the BD triad $A, \, A', \, A''$ is thin.  Let $B$ denote the linear transformation from Lemma \ref{defB}.  Let $B', \, B''$ denote the linear transformations from Definition \ref{defB'B''}.  With reference to the irreducible $\boxtimes$-module structure on $V$ constructed in the proof of Theorem \ref{thm:BDtriadtet} we make the following observations.  By Theorem \ref{thm:tetBDtriad} each of the following 
\begin{eqnarray*}
A, \, A', \, A'' \qquad \qquad -A'', \, B', \, -B \qquad \qquad B, \, -A', \, -B'' \qquad \qquad B'', \, -B', \, -A  
\end{eqnarray*}
acts on $V$ as a reduced BD triad.  By \cite[Theorem 4.8]{Funk-Neubauer17}, \cite[Theorem 3.8]{Hartwig07}, and \cite[Corollary 12.4]{HarTer07} each of the following
\begin{eqnarray*}
A', \, -A'', \, B \qquad \qquad A, \, -A', \, B'' \qquad \qquad A, \, -A'', \, -B' \qquad \qquad B', \, B, \, B''  
\end{eqnarray*}
acts on $V$ as a reduced BD triple of base $1$.  From this and \cite[Lemma 5.8]{Funk-Neubauer17} each of the following  
\begin{eqnarray*}
A', \, -A'' \qquad \qquad -A'', \, B \qquad \qquad B, \, A' \qquad \qquad A, \, -A'  \qquad \qquad -A', \, B'' \qquad \qquad B'', \, A  \\
A, \, -A'' \qquad \qquad -A'', \, -B' \qquad \qquad -B', \, A \qquad \qquad B', \, B  \qquad \qquad B, \, B'' \qquad \qquad B'', \, B' 
\end{eqnarray*}
acts on $V$ as a reduced BD pair of base $1$.  By \cite[Theorems 1.7, 2.4]{Hartwig07} each of the following 
\begin{eqnarray*}
A, \, B \qquad \qquad A', \, B' \qquad \qquad A'', \, B''  
\end{eqnarray*}
acts on $V$ as a tridiagonal pair of Krawtchouk type \cite{ItoTer074}.
\end{note}

\section{A Counterexample}

In this section we provide a counterexample to show that the thin assumption in Theorem \ref{thm:BDtriadtet} is needed.  That is, we demonstrate that the conclusion of Theorem \ref{thm:BDtriadtet} can fail to hold for BD triads which are not thin.  Consider the vector space $\mathbb{K}^6$ and let $A: \mathbb{K}^6 \to \mathbb{K}^6$, $A' :\mathbb{K}^6 \to \mathbb{K}^6$, and $A'' :\mathbb{K}^6 \to \mathbb{K}^6$ denote the linear transformations given by the following matrices:
\begin{eqnarray*}
A = \begin{bmatrix}
-3 & 0 & 0 & 0 & 0 & 0 \\ 1 & -1 & 0 & 0 & 0 & 0 \\ 1 & 0 & -1 & 0 & 0 & 0 \\ 0 & 2 & 0 & 1 & 0 & 0 \\ 0 & 1 & 2 & 0 & 1 & 0 \\ 0 & 0 & 0 & 3 & 0 & 3
\end{bmatrix}
A' = \begin{bmatrix}
-3 & 0 & 0 & 0 & 0 & 0 \\ -2 & -1 & 0 & 0 & 0 & 0 \\ 0 & 0 & -1 & 0 & 0 & 0 \\ 0 & -4 & 0 & 1 & 0 & 0 \\ 0 & 0 & -2 & 0 & 1 & 0 \\ 0 & 0 & 0 & -6 & 0 & 3
\end{bmatrix}
A'' = \begin{bmatrix}
-3 & 0 & 0 & 0 & 0 & 0 \\ 0 & -1 & 0 & 0 & 0 & 0 \\ 0 & 0 & -1 & 0 & 0 & 0 \\ 0 & 0 & 0 & 1 & 0 & 0 \\ 0 & 0 & 0 & 0 & 1 & 0 \\ 0 & 0 & 0 & 0 & 0 & 3
\end{bmatrix}.
\end{eqnarray*}
It is straightforward to verify that $A, \, A', \, A''$ is a reduced BD triad on  $\mathbb{K}^6$ whose first, second, and third eigenvalue sequences are each $\lbrace 2i-3 \rbrace _{i=0}^3$.  Notice that the shape of $A, \, A', \, A''$ is $(1,2,2,1)$, and so this BD triad is not thin.  Consider the corner triad $X_{03}, X_{13}, X_{23}$ in $\boxtimes$.  Suppose, towards a contradiction, that the conclusion of Theorem \ref{thm:BDtriadtet} holds.  That is, suppose there exists a $\boxtimes$-module structure on $\mathbb{K}^6$ such that $(X_{03} - A)\mathbb{K}^6=0$, $(X_{13} - A')\mathbb{K}^6=0$, and $(X_{23} - A'')\mathbb{K}^6=0$.  Then the matrix representing $X_{02}$ is given by
\begin{eqnarray*}
X_{02} = \begin{bmatrix}
3 & 12 & 0 & 0 & 0 & 0 \\ 0 & 1 & 0 & 8 & 0 & 0 \\ 0 & 0 & 1 & 5 & 2 & 0 \\ 0 & 0 & 0 & -1 & 0 & 4 \\ 0 & 0 & 0 & 0 & -1 & 6 \\ 0 & 0 & 0 & 0 & 0 & -3
\end{bmatrix}.
\end{eqnarray*}
However, direct calculation shows that 
$$[X_{13}, [X_{13}, [X_{13}, X_{02}]]] \neq 4[X_{13}, X_{02}].$$
which contradicts Definition \ref{def:tet}(iii).  Thus, there does not exist a $\boxtimes$-module structure on $\mathbb{K}^6$ such that $(X_{03} - A)\mathbb{K}^6=0$, $(X_{13} - A')\mathbb{K}^6=0$, and $(X_{23} - A'')\mathbb{K}^6=0$, and so the conclusion of Theorem \ref{thm:BDtriadtet} does not hold. 

\section{Acknowledgment}
The author would like to thank Professor Paul Terwilliger for his advice and many helpful suggestions.  His comments greatly improved this paper.  In particular, his comment indicating the relevance of \cite{ItoTerinpress4} to this paper was especially helpful.  \\

\bibliographystyle{amsplain}
\bibliography{FunkNeubauer}

\noindent
Darren Funk-Neubauer \hfil\break
\noindent Department of Mathematics and Physics \hfil\break
\noindent Colorado State University - Pueblo \hfil\break
\noindent 2200 Bonforte Boulevard \hfil\break
\noindent Pueblo, CO 81001 USA \hfil\break
email:   {\tt darren.funkneubauer@csupueblo.edu} \hfil\break
phone:  (719) 549 - 2693 \hfil\break
fax:  (719) 549 - 2962

\end{document}